\documentclass{imsart}

\RequirePackage{amsthm,amsmath,amsfonts,amssymb}
\RequirePackage[numbers,sort&compress]{natbib}
\RequirePackage[colorlinks,citecolor=blue,urlcolor=blue]{hyperref}
\usepackage{mathrsfs}
\usepackage{dsfont}
\usepackage{lmodern}
\usepackage{color}
\usepackage{bbm}
\usepackage[latin1]{inputenc}
\usepackage{csquotes}
\usepackage{graphicx,epstopdf}
\usepackage[dvips]{epsfig}
\usepackage{latexsym}
\usepackage{exscale}
\usepackage[toc]{appendix}
\usepackage{diagbox}
\usepackage{multirow}
\usepackage{tcolorbox}
\usepackage{xcolor}
\usepackage{soul}
\usepackage{ulem}
\usepackage{booktabs}
\usepackage[ruled]{algorithm2e}

\usepackage{tabularx}
\usepackage{array}
\usepackage{colortbl}
\usepackage{float}
\startlocaldefs
\newtheorem{theorem}{Theorem}[section]

\newtheorem{alg}[theorem]{Algorithm}

\theoremstyle{remark}

\newtheorem{definition}[theorem]{Definition}

\DeclareMathOperator*{\argmin}{arg\,min}

\numberwithin{equation}{section}
\endlocaldefs
\begin{document}
\begin{frontmatter}
\title{Estimation of the Self-similarity Index of Non-stationary Increments Self-similar Processes via Lamperti Transformations}
\runtitle{Estimation of the Self-similarity Index}

\begin{aug}
\author[A]{\fnms{William} \snm{Wu}
\ead[label=e1]{wwu49387@gmail.com}}
\and
\author[B]{\fnms{Qidi} \snm{Peng}\ead[label=e2,mark]{qidi.peng@cgu.edu}}


\address[A]{Rancho Cucamonga High School,
\printead{e1}}

\address[B]{Institute of Mathematical Sciences, Claremont Graduate University,
\printead{e2}}

\runauthor{W. Wu and Q. Peng}
\end{aug}

\begin{abstract}
  We introduce a novel method for estimating the self-similarity index of a general $H$-self-similar process with either stationary or non-stationary increments. The estimation algorithm is developed based on a modified Lamperti transformation, which transforms $H$-self-similar processes to stationary ones. As an application, we show how to use this approach to estimate the self-similarity index of fractional Brownian motion, subfractional Brownian motion, bifractional Brownian motion, and trifractional Brownian motion. Simulation study is performed to support the consistency of our estimators. Implementation in Python is publicly shared. Application on the estimation of the self-similarity index of the Nile river water level data from the year 900 to 1200 C.E..
\end{abstract}


\end{frontmatter}


\section{Literature Review on the Estimation of the Self-similarity Index}
\label{sec:introduction:selfsimilar}
Self-similar processes are stochastic processes that exhibit the property of invariance in distribution under scaling of time and space. These processes are essential notions in stochastic processes and fractal geometry \citep{Samorodnitsky1994}, with diverse applications in real-world fields such as packet inter-arrival times and burst lengths \citep{leland1994self,abry2000self}, and global financial equity markets \citep{Bianchi2008,Bianchi2012,peng2011inference,Peng2018}. Such processes exhibit the phenomenon of \enquote{self-similarity} and are often associated with the ``rough path" feature: the process's sample path is mainly characterized by its self-similarity index. In this paper, we consider the $H$-self-similar process defined below:
\begin{definition}
\label{def:selfsimilar}
We say $V=\{V(t)\}_{t\in \mathbb{R}}$ is an $H$-self-similar process with self-similarity index $H\in(0,1)$, if
\begin{equation}
\label{self-similar}
    \left\{V(at)\right\}_{t\in\mathbb R}\stackrel{\mbox{f.d.d.}}{=}\left\{|a|^HV(t)\right\}_{t\in\mathbb R}, ~\mbox{for any real number}~ a\ne0,
\end{equation}
where $\stackrel{\mbox{f.d.d.}}{=}$ denotes equality in finite dimensional distribution, i.e., for any $n\ge1$ and any $t_1,\ldots,t_n\in\mathbb R$, 
$$
    \left(V(at_1),\ldots,V(at_n)\right)\stackrel{\mbox{law}}{=}|a|^H\left(V(t_1),\ldots,V(t_n)\right), ~\mbox{for any real number}~ a\ne0.
$$
\end{definition}

 Concerning to estimate the self-similarity index $H$ of a $H$-self-similar process, there exist several method in the literature, however, each of them deals with a special case. Therefore, the goal of this paper is to remedy this inconvenience by developing a general method to estimate the self-similarity index of a large class of fourth order $H$-self-similar processes. 

\section{A Modified Lamperti Transformation: A Novel Estimation Approach}
\subsection{Lamperti Transformation}
Our simulation approaches are inspired by the Lamperti transformations. In 1962, Lamperti introduced the Lamperti transformation to transform $H$-self-similar processes to stationary processes, and the inverse operator does so vice versa. The two transformations are defined as follows \citep{lamperti1962semi}.
\begin{definition}
\label{def:Lamperti}
Let $\{V(t)\}_{t\in \mathbb{R}}$ be an $H$-self-similar process with self-similarity index $H\in(0,1)$, the Lamperti transformation of $\{V(t)\}_{t\in \mathbb{R}}$ is defined to be
\begin{equation}
\label{Lamperti}
    \mathcal L_V(t) = e^{-tH}V(e^{t}) ~\mbox{for}~ t\in\mathbb R.
\end{equation}
Let $\{U(t)\}_{t\in \mathbb{R}}$ be a stationary stochastic process, i.e., for any integer $n\ge1$, any time indices $t_1,\ldots,t_n\in\mathbb R$ and any time shifting parameter $h>0$, we have
$$
\left(U_{t_1},\ldots,U_{t_n}\right) \stackrel{\mbox{law}}{=}\left(U_{t_1+h},\ldots,U_{t_n+h}\right).
$$
For each $H\in(0,1)$, the inverse Lamperti transformation of $\{U(t)\}_{t\in \mathbb{R}}$ is given by
$$
\mathcal L_U^{-1}(0) = 0;~\mathcal L_U^{-1}(t) = t^H U(\log t)~\mbox{for}~ t>0.
$$
\end{definition}
Through the paper, we assume that 
\begin{equation}
\label{variance}
\sigma^2=\mathbb Var\left(V(1)\right).
\end{equation}
The main goal of this paper is to provide a novel and general method to estimate the self-similarity index of $\{V(t)\}_{t\in[0,1]}$. The options of the methodologies depend on the situation for the cases for known $\sigma$ and unknown $\sigma$. 
In Definition \ref{def:Lamperti}, the fact that $\{\mathcal L_V(t)\}_{t\in\mathbb R}$ is stationary inspires us to develop a novel estimation strategy. Heuristically speaking, our method involves first transforming the $H$-self-similar process to an ergodic stationary one by using some modified Lamperti transformation, then running statistical estimation of the self-similarity index on the latter process.
\subsection{Estimation of the Self-similarity Index When $\sigma^2$ is Given}
If $\sigma^2$ is known, we assume that $\{\mathcal L_V(t)\}_{t\in[0,1]}$ is ergodic. In view of the definition of $H$-self-similar process and the Lamperti transformation, we suggest the following way to estimate the self-similarity index $H$ of the process $\{V(t)\}_{t\in[0,1]}$. 
 \begin{alg}\ \\
 \label{alg:WP_known_sigma}
 \begin{description}
 \item[Step 1] Define the sequences 
  $a=(a_0,\ldots,a_n)
 $ and $b=(b_0,\ldots,b_n)$: 
\begin{equation}
\label{def:a}
a=\left(V\left(\frac{\lfloor n^{j/n}\rfloor}{n}\right)\right)_{j=0,\ldots,n}~\mbox{and}~b=\left(n^{1-j/n}\right)_{j=0,\ldots,n},
\end{equation}
where  $\lfloor\bullet\rfloor$ denotes the floor number. 
  \item[Step 2] Use a non-linear least squares method to solve $H$ from the equation
\begin{equation}
\label{H:equation}
f_n(H):=\frac{\sum_{j=0}^na_j^2b_j^{2H}}{n+1}-\sigma^2=0.
\end{equation}
For example, the solution $\widehat H$ can be solved from Halley's method \cite{scavo1995geometry}.
\end{description}
\end{alg}
Let us heuristically explain the rationales of Algorithm \ref{alg:WP_known_sigma}.
\begin{enumerate}
\item In \textbf{Step 1},  it is from the integer $J_n = \lceil-n\log_n(n^{1/n}-1)\rceil
 $ such that $\{\lfloor n^{j/n}\rfloor\}_{j\in\{0,\ldots,n\}}$ are all distinct. In fact, $J_n$ is the smallest integer 
such that
$$
n^{(j+1)/n}- n^{j/n}\ge1~\mbox{for all}~ j\ge J_n.
$$
\item In \textbf{Step 2}, we have used the fact that
\begin{equation}
\label{def:inverse_Lamperti}
U_{n}\left(\frac{j}{n}\right):=\mathcal L_V\left(\left(\frac{j}{n}-1\right)\log (n)\right)=n^{-H(j/n-1)}V(n^{j/n-1})~\mbox{for}~j\in\left\{0,\ldots,n\right\}
\end{equation}
is stationary and ergodic. Therefore, 
\begin{equation}
\label{convergence_Proba}
\frac{\sum_{j=0}^{n}V(n^{j/n-1})^2n^{2H(1-j/n)}}{n+1}-\mathbb Var(V(1))\xrightarrow[n\to+\infty]{\mathbb P}0.
\end{equation}
We also have
$$
V\left(n^{j/n-1}\right)\approx V\left(\frac{\lfloor n^{j/n}\rfloor}{n}\right)=a_j
~\mbox{a.s..}
$$
It results that
\begin{equation}
\label{convergence_Proba_estimation}
\frac{\sum_{j=0}^{n}a_j^2b_j^{2H}}{n+1}-\sigma^2\xrightarrow[n\to+\infty]{\mathbb P}0.
\end{equation} 
\end{enumerate}
\subsection{Estimation of the Self-similarity Index When $\sigma^2$ is Unknown}
If $\sigma^2$ is unknown, we further assume that $\mathbb E(V(1)^4)<+\infty$. We then observe that
$$
\widehat\kappa_n(H):=\frac{(n+1)\sum_{j=0}^{n}a_j^4b_j^{4H}}{\left(\sum_{j=0}^{n}a_j^2b_j^{2H}\right)^2}\xrightarrow[n\to+\infty]{\mbox{a.s.}}\frac{\mathbb E (V(1)^4)}{(\mathbb E(V(1)^2))^2}.
$$
$\widehat\kappa_n(H)$ is then a consistent estimator of the 
kurtosis of the stationary distribution $V(1)$. This kurtosis is the measure of dispersion of $Z^2$ around its expectation $1$, where $Z=V(1)^2/\sqrt{\mathbb E(V(1)^2)}$. Therefore, the maximum likelihood estimation method is equivalent to minimizing the kurtosis via $H$, i.e.,
$$
\widehat H_n:=\argmin_{H\in(0,1)}\widehat\kappa_n(H)\xrightarrow[n\to+\infty]{\mbox{a.s.}}H.
$$
We then propose the following algorithm:
 \begin{alg}\ \\
 \label{alg:WP_unknown_sigma}
 \begin{description}
 \item[Step 1] Define the sequences 
  $a=(a_0,\ldots,a_n)
 $ and $b=(b_0,\ldots,b_n)$ as in  (\ref{def:a}).  
  \item[Step 2] Use Brent's method \cite{brent1971algorithm} to solve the following optimization problems
\begin{equation}
\label{H:equation_2}
\widehat H_n=\argmin\limits_{H\in(0,1)}\frac{\sum_{j=0}^{n}a_j^4b_j^{4H}}{\left(\sum_{j=0}^{n}a_j^2b_j^{2H}\right)^2}.
\end{equation}
\end{description}
\end{alg}

\section{Estimation of the Self-similarity Index with Examples}
In this section, we will estimate and the self-similarity index for fractional Brownian motion, subfractional Brownian motion, bifractional Brownian motion, and trifractional Brownian motion. 
\subsection{Fractional Brownian motion}
     The fractional Brownian motion $\{B^H(t)\}_{t\ge0}$ \cite{Mandelbrot1968} with the Hurst parameter $H\in(0,1)$ that we consider in this framework is defined as a zero mean Gaussian process with $B^H(0)=0$ a.s. and covariance function given by: for $s,t\ge0$,
\begin{equation}
\label{cov_fBm}
\gamma_{B^H}(s,t)=\frac{\sigma^2}{2}\left(|t|^{2H}+|s|^{2H}-|t-s|^{2H}\right).
\end{equation}
In a simulation study, we use our estimation for processes simulated using Wood Chan's method \cite{WoodChan}. Simulations and estimations were performed. We generate $200$ trajectories of fBm with length $N=2^l$ for $l\in\{7,8,9,10\}$ and $H\in\{0.2,0.5,0.7,0.8\}$. All of the following results show the mean estimated $H$ and the mean squared error in parentheses. S

\begin{table}[h]
\centering
\caption{Estimator of $H$ with MSE based on fBm simulated by Wood-Chan's Method: $\sigma^2=1$}
\label{tab:wc-fbm-halley}
\begin{tabular}{c cccc}
\toprule
$H$ & $N=128$ & $N=256$ & $N=512$ & $N=1024$ \\
\midrule
0.2 & 0.237 (0.0213) & 0.226 (0.0191) & 0.220 (0.0121) & 0.209 (0.0084) \\
0.5 & 0.538 (0.0524) & 0.505 (0.0195) & 0.508 (0.0149) & 0.510 (0.0112) \\
0.7 & 0.733 (0.0715) & 0.686 (0.0502) & 0.715 (0.0211) & 0.715 (0.0125) \\
0.8 & 0.823 (0.0916) & 0.849 (0.0361) & 0.809 (0.0238) & 0.818 (0.0219) \\
\bottomrule
\end{tabular}
\end{table}

The following results are for Algorithm \ref{alg:WP_unknown_sigma} where the self-similar scaling parameter is not known. 

\begin{table}[H]
\centering
\caption{Estimator of $H$ with MSE based on fBm simulated by Wood-Chan's Method: $\sigma^2$ is unknown}
\label{tab:dpw-fbm-convergence}
\begin{tabular}{c ccccc}
\toprule
$H$ & $N=128$ & $N=256$ & $N=512$ & $N=1024$ & $N=8192$ \\
\midrule
0.2 &
0.286 (0.0405) &
0.264 (0.0307) &
0.235 (0.0205) &
0.228 (0.0112) &
0.219 (0.0067) \\
0.5 &
0.504 (0.0503) &
0.504 (0.0311) &
0.496 (0.0234) &
0.515 (0.0203) &
0.507 (0.0105) \\
0.7 &
0.668 (0.0350) &
0.684 (0.0317) &
0.683 (0.0255) &
0.658 (0.0253) &
0.691 (0.0113) \\
0.8 &
0.743 (0.0401) &
0.721 (0.0466) &
0.754 (0.0276) &
0.753 (0.0245) &
0.773 (0.0148) \\
\bottomrule
\end{tabular}
\end{table}

For comparison, we can test the accuracy against the quadratic variation estimation method (QV estimator) \cite{Coeurjolly2005} and the complex variation estimation method \cite{istas2012estimating} by running the same simulations. The results are given below.

\begin{table}[h]
\centering
\caption{Quadratic Variation Estimator of $H$.}
\label{tab:qv-wc-fbm}
\begin{tabular}{c cccc}
\toprule
$H$ & $N=128$ & $N=256$ & $N=512$ & $N=1024$ \\
\midrule
0.2 & 0.193 (0.0056) & 0.204 (0.0025) & 0.207 (0.0012) & 0.199 (0.0006) \\
0.5 & 0.499 (0.0072) & 0.498 (0.0030) & 0.499 (0.0015) & 0.497 (0.0009) \\
0.7 & 0.683 (0.0080) & 0.690 (0.0036) & 0.702 (0.0021) & 0.698 (0.0009) \\
0.8 & 0.794 (0.0088) & 0.799 (0.0044) & 0.797 (0.0020) & 0.794 (0.0010) \\
\bottomrule
\end{tabular}
\end{table}

\begin{table}[h]
\centering
\caption{Complex Variations Estimator of $H$.}
\label{tab:cv-wc-fbm}
\begin{tabular}{c ccc}
\toprule
$H$ & $N=128$ & $N=256$ & $N=512$ \\
\midrule
0.5 & 0.500 (0.050) & 0.500 (0.022) & 0.490 (0.012) \\
0.7 & 0.680 (0.040) & 0.690 (0.023) & 0.710 (0.010) \\
0.8 & 0.825 (0.037) & 0.805 (0.020) & 0.800 (0.008) \\
\bottomrule
\end{tabular}
\end{table}
The comparison shows that our estimation methods are competitive to the literature ones for estimating the Hurst parameter of fractional Brownian motion.

Additional visualization of the error analysis result is given below.
\newpage
\begin{figure}[H]
\centering
\includegraphics[width=0.45\linewidth]{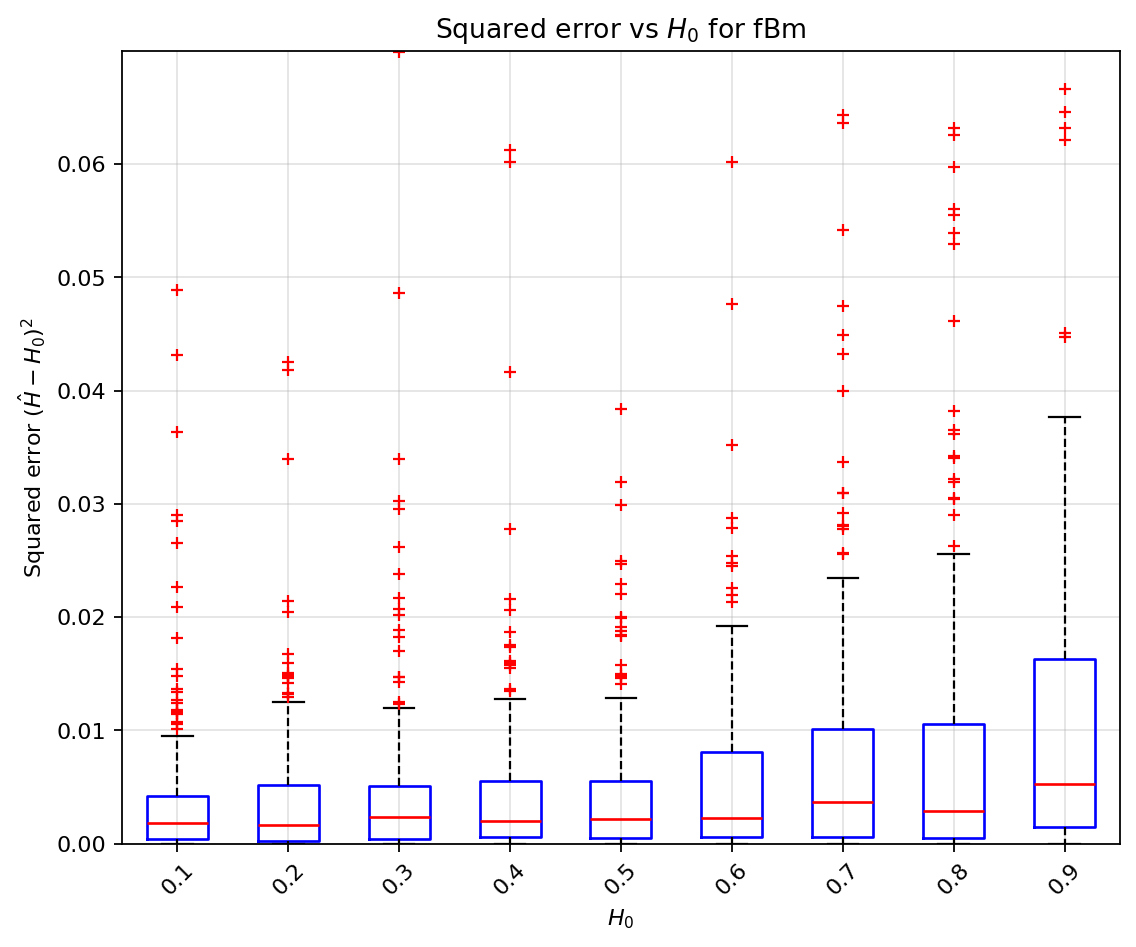}
\includegraphics[width=0.45\linewidth]{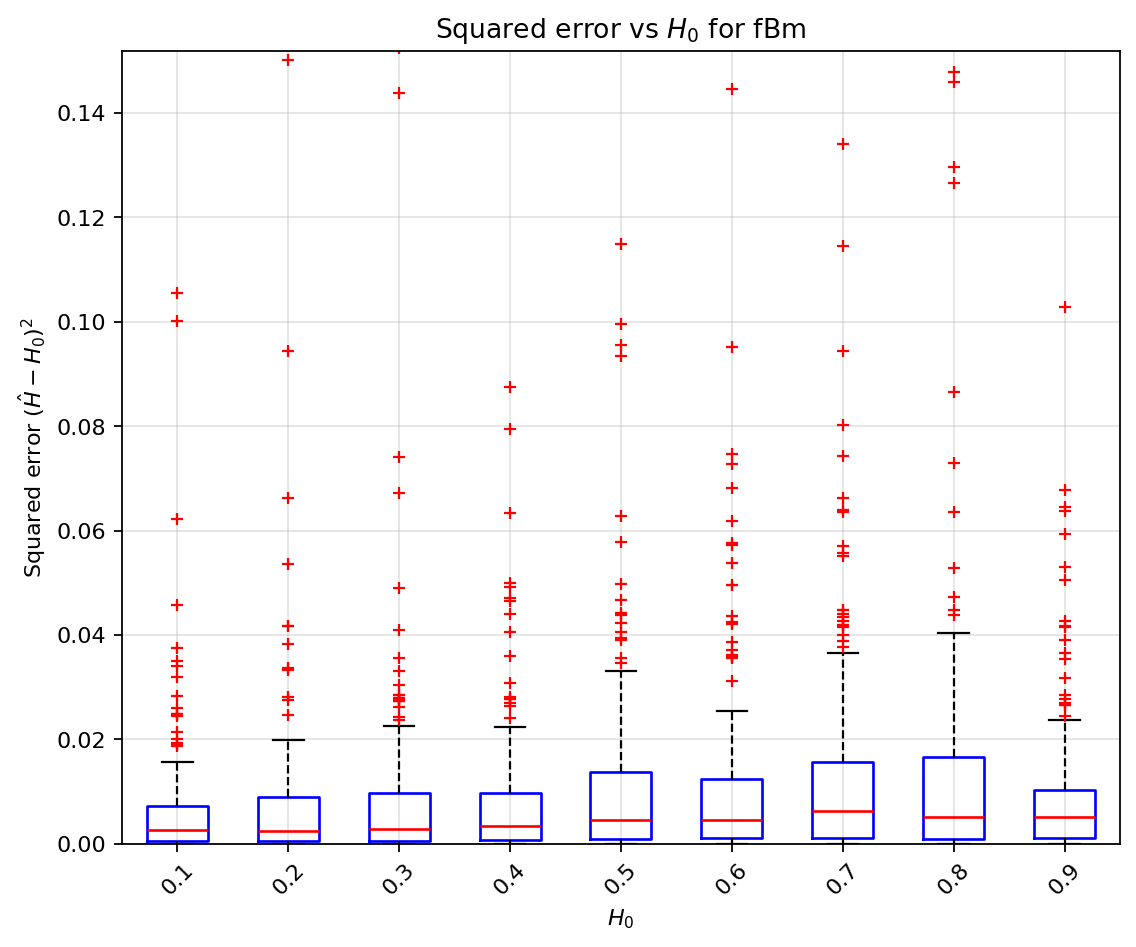}
\caption{The results are given by left: Algorithm \ref{alg:WP_known_sigma}, right: Algorithm \ref{alg:WP_unknown_sigma}. Each $H_0$ consists of 200 simulated paths of length $N=4096$. Squared error of the estimated parameter is computed for each path and drawn into a boxplot. Outliers over 99th percentile are cut out for the sake of clear visuals.}
\end{figure}
\begin{figure}[H]
\centering
\includegraphics[width=0.8\linewidth]{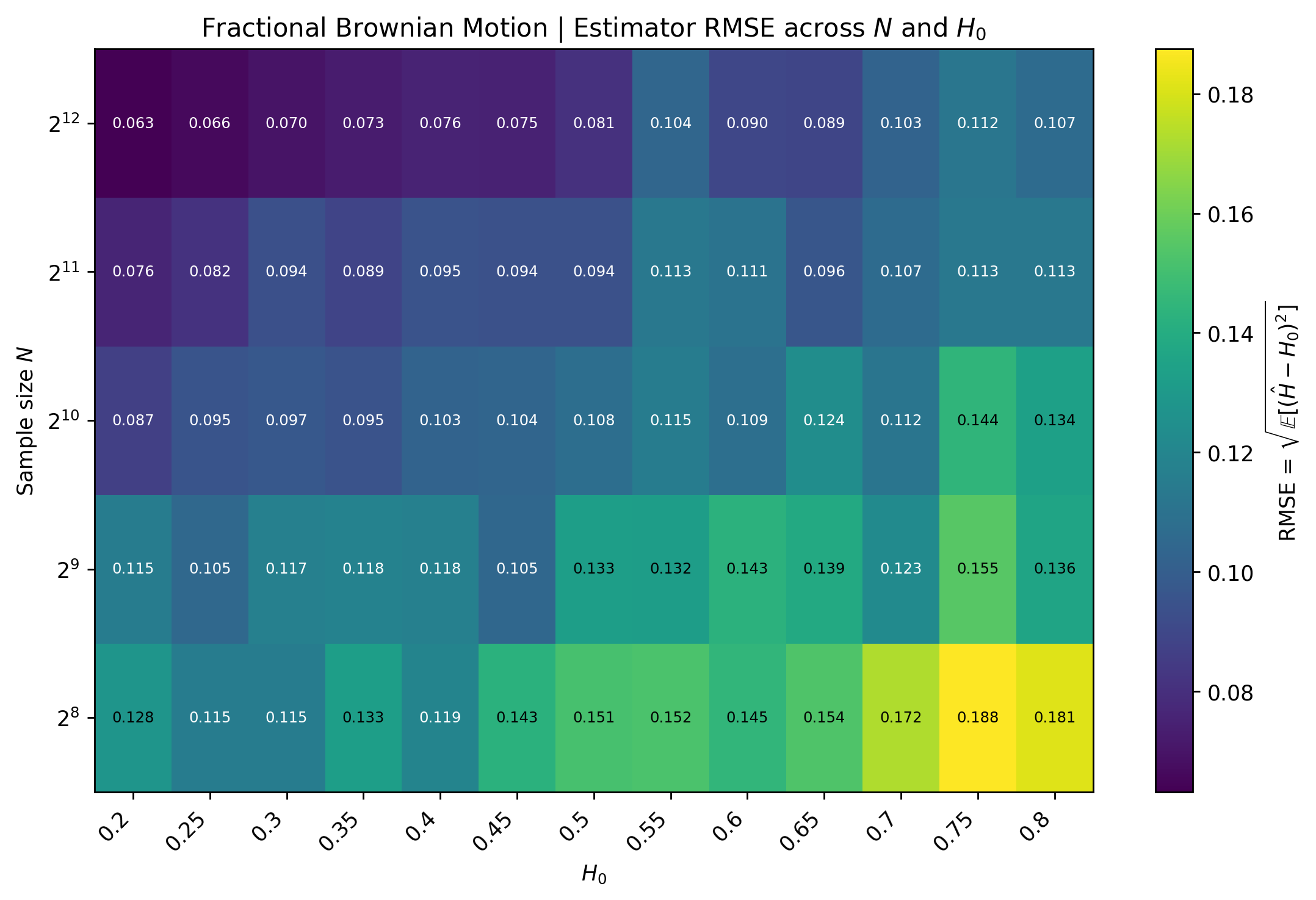}
\caption{200 paths of fBm are simulated for each corresponding $H_0$ and $N$ value with the RMSE plotted as a heatmap distribution. The results show empirical convergence towards 0 as the path length increases and $H_0$ value decreases.}
\label{fig:boxplot_fbm}
\end{figure}
The above figures show that in both cases, the overall performance of the estimation is consistent.

\subsection{Subfractional Brownian motion}
The subfractional Brownian motion $\{S^H(t)\}_{t\ge0}$ is defined to be a
zero-mean Gaussian process with $S^H(0)=0$ a.s. and covariance function (see \cite{bojdecki2004sub}): for $s,t\ge0$,
\begin{equation}
\label{def:sfBm}
\mathbb Cov\left(S^H(t),S^H(s)\right)=\sigma^2\left(t^{2H}+s^{2H}-\frac{1}{2}\left\{(t+s)^{2H}+|t-s|^{2H}\right\}\right),
\end{equation}
where the parameter $H\in(0,1)$. By (\ref{def:sfBm}), we know that $\{S^H\}_{t\ge0}$ is an $H$-self-similar process. However, unlike the fractional Brownian motion, the subfractional Brownian motion has non-stationary increments when $H\ne 1/2$. This property makes its simulation quite challenging. The simulation approach existing in the literature is based on the DPW method provided in \cite{peng2025investigation} with available code in  \url{https://pypi.org/project/fractal-analysis/}. 

       To estimate $H$ when $\sigma^2$ is known as $1$, it suffices to note that
       $$
       \sigma^2=Var(S^H(1))=2-2^{2H-1}.
       $$
       Modifying a bit Algorithm \ref{alg:WP_known_sigma}, we obtain the algorithm to estimate $H$.
\begin{alg}
 \label{alg:sfBm}
Let $(S^H(j/n))_{j=1,\ldots,n}$ be an observed sample path of $\{S^H(t)\}_{t\in[0,1]}$ with covariance function $\gamma_{S^H}$ with $\sigma^2=1$. 
 \begin{description}

 \item[Step 1] Let the sequences 
  $a=(a_0,\ldots,a_n)
 $ and $b=(b_0,\ldots,b_n)$: 
$$
a=\left(S^H\left(\frac{\lfloor n^{j/n}\rfloor}{n}\right)\right)_{j=0,\ldots,n}~\mbox{and}~b=\left(n^{1-j/n}\right)_{j=0,\ldots,n},
$$
where  $\lfloor\bullet\rfloor$ denotes the floor number. 
  \item[Step 2] Use a non-linear least squares method to solve $H$ from the equations
$$
f_n(H):=\frac{\sum_{j=0}^{n}a_j^2b_j^{2H}}{n+1}-(2-2^{2H-1})=0.
$$
\end{description}
\end{alg}
We show that our estimation methods perform well for sfBm using the DPW simulation algorithm. We generate 200 trajectories of sfBm with length $N=2^l$ for $l\in\{7,8,9,10\}$ and $H\in\{0.2,0.5,0.7,0.8\}$. The following results show the mean estimator $H$ with the MSE in parentheses. 

\begin{table}[H]
\centering
\caption{Subfractional Brownian motion: Algorithm \ref{alg:sfBm}}
\label{tab:sfbm-halley}
\begin{tabular}{c cccc}
\toprule
$H$ & $N=128$ & $N=256$ & $N=512$ & $N=1024$ \\
\midrule
0.2 & 0.228 (0.0218) & 0.216 (0.0128) & 0.228 (0.0090) & 0.214 (0.0071) \\
0.5 & 0.509 (0.0142) & 0.509 (0.0107) & 0.508 (0.0077) & 0.508 (0.0060) \\
0.7 & 0.695 (0.0081) & 0.703 (0.0066) & 0.702 (0.0052) & 0.705 (0.0044) \\
0.8 & 0.785 (0.0061) & 0.799 (0.0051) & 0.799 (0.0033) & 0.801 (0.0027) \\
\bottomrule
\end{tabular}
\end{table}
The simulation results for algorithm \ref{alg:WP_unknown_sigma} is given below:

\begin{table}[H]
\centering
\caption{Subfractional Brownian motion: Algorithm \ref{alg:WP_unknown_sigma}}
\label{tab:dpw-sfbm-convergence}
\begin{tabular}{c ccccc}
\toprule
$H$ & $N=128$ & $N=256$ & $N=512$ & $N=1024$ & $N=8192$ \\
\midrule
0.2 &
0.251 (0.0370) &
0.266 (0.0272) &
0.250 (0.0274) &
0.247 (0.0191) &
0.220 (0.0086) \\
0.5 &
0.500 (0.0388) &
0.501 (0.0293) &
0.516 (0.0227) &
0.500 (0.0153) &
0.500 (0.0096) \\
0.7 &
0.635 (0.0494) &
0.657 (0.0273) &
0.661 (0.0263) &
0.668 (0.0240) &
0.680 (0.0109) \\
0.8 &
0.750 (0.0351) &
0.761 (0.0314) &
0.762 (0.0268) &
0.732 (0.0291) &
0.762 (0.0171) \\
\bottomrule
\end{tabular}
\end{table}

Additional data visualization is given. 

\begin{figure}[H]
\centering
\includegraphics[width=0.45\linewidth]{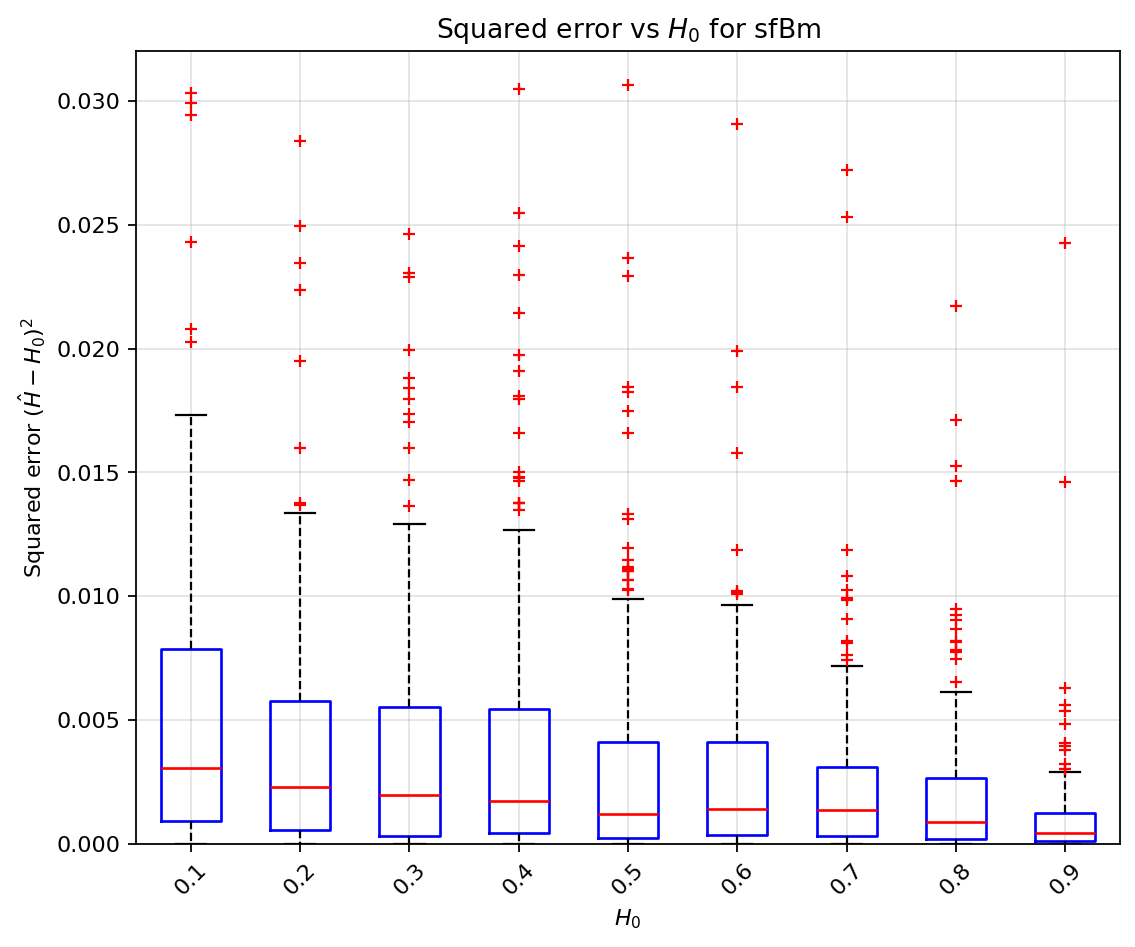}
\includegraphics[width=0.45\linewidth]{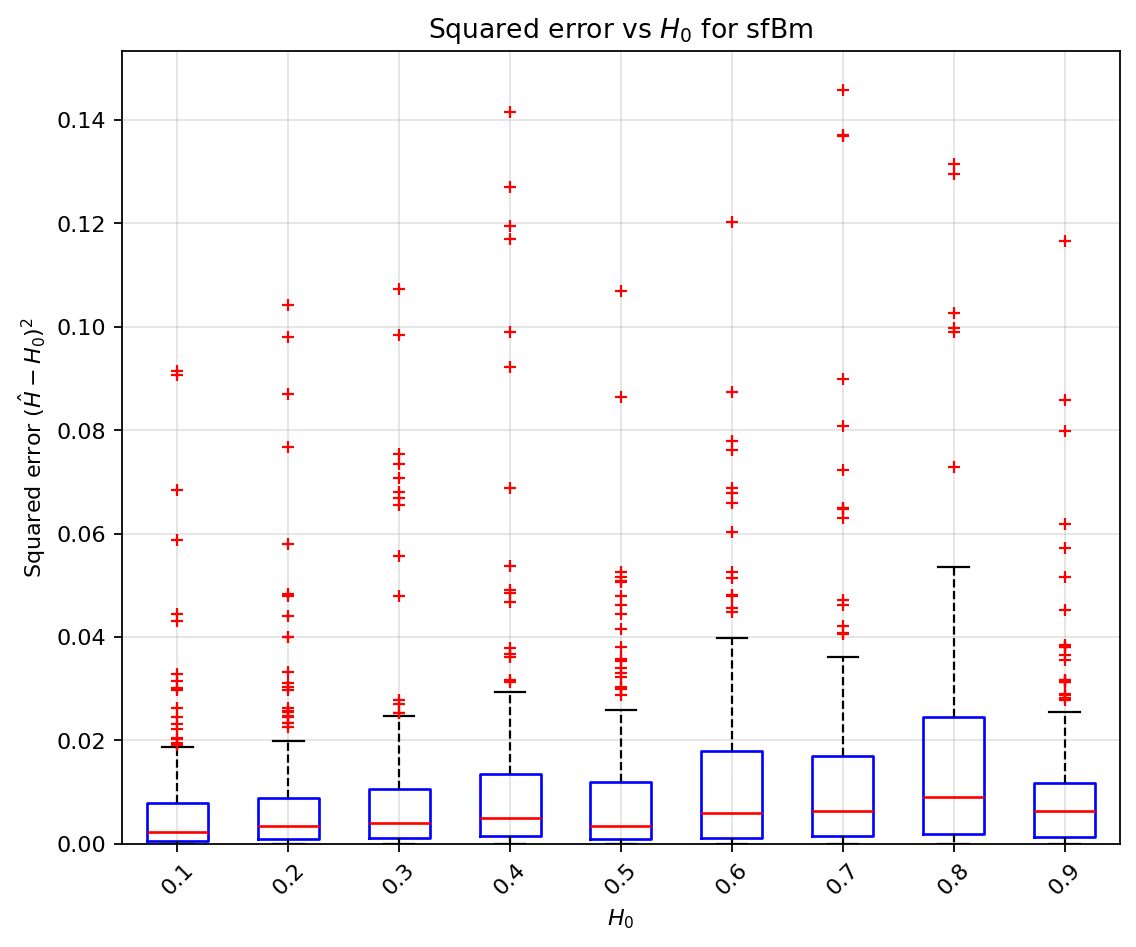}
\caption{The results are given by left: Algorithm \ref{alg:WP_known_sigma}, right: Algorithm \ref{alg:WP_unknown_sigma}. Each $H_0$ consists of 200 simulated paths of length $N=4096$. Squared error of the estimated parameter is computed for each path and drawn into a boxplot. Outliers over 99th percentile are cut out for the sake of clear visuals.}
\label{fig:boxplot_sfbm_1}
\end{figure}

\begin{figure}[H]
\centering
\includegraphics[width=0.8\linewidth]{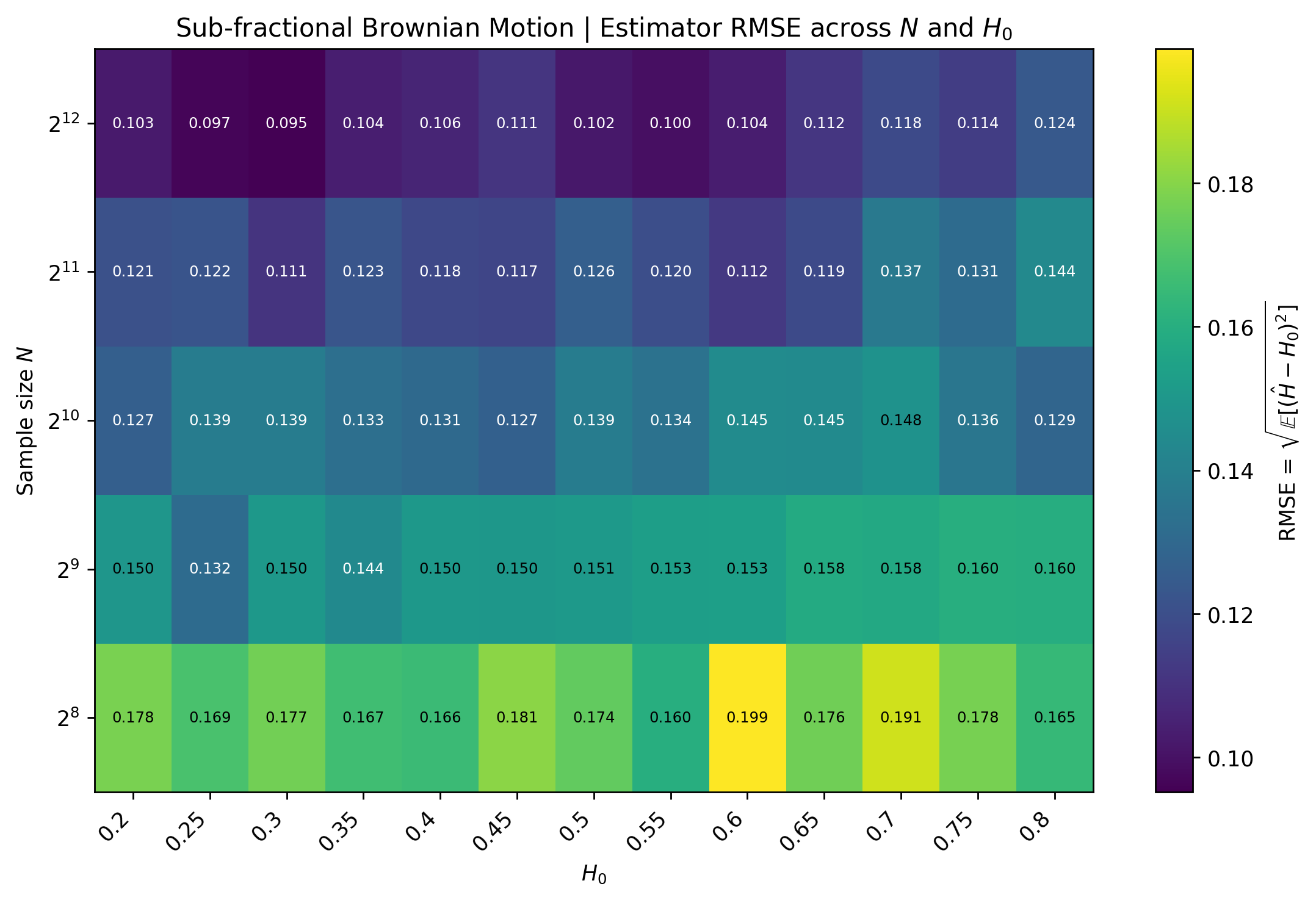}
\caption{200 paths of sfBm are simulated for each corresponding $H_0$ and $N$ value with the RMSE plotted as a heatmap distribution. The results show empirical convergence towards 0 as the path length increases.}
\label{fig:boxplot_sfbm_2}
\end{figure}

\subsection{Bifractional Brownian motion}
The bifractional Brownian motion $\{B^{H,K}(t)\}_{t\ge0}$ is a
zero-mean Gaussian process with $B^{H,K}(0)=0$ a.s. and covariance function (see \cite{houdre2003example}): for $s,t\ge0$,
\begin{equation}
\label{def:bfBm}
\mathbb Cov\left(B^{H,K}(t),B^{H,K}(s)\right)=\frac{\sigma^2}{2^K}\left\{(t^{2H}+s^{2H})^K-|t-s|^{2HK}\right\},
\end{equation}
where the parameters $H\in(0,1)$ and $K\in(0,1]$. By (\ref{def:bfBm}), we see that $\{B^{H,K}(t)\}_{t\ge0}$ is an $HK$-self-similar process. When $K=1$, $\{B^{H,1}(t)\}_{t\ge0}$ becomes a fractional Brownian motion with the Hurst parameter $H$. For $K<1$, the bifractional Brownian motion does not have stationary increments.  To estimate the self-similarity $HK$, it suffices to apply Algorithm \ref{alg:WP_known_sigma} and Algorithm \ref{alg:WP_unknown_sigma}. Note that the simulation method is also provided in the DPW method \cite{peng2025investigation}. 

The following table shows our estimation method for bfBm using the DPW simulation algorithm. Like before, we will generate $200$ trajectories of bfBm with lengths $N=128,256,512,1024$. This will be tested on paths with $(H,K)\in(\{0.2,0.8\},\{0.5,0.8\})$ to estimate the self-similarity value $HK$. 

\begin{table}[H]
\centering
\caption{Bifractional Brownian motion: Algorithm \ref{alg:WP_known_sigma} with $\sigma^2=1$}
\label{tab:bfbm}
\begin{tabular}{c cccc}
\toprule
$(H,K)$ & $N=128$ & $N=256$ & $N=512$ & $N=1024$ \\
\midrule
(0.2, 0.5) & 0.108 (0.0074) & 0.107 (0.0058) & 0.114 (0.0044) & 0.111 (0.0026) \\
(0.2, 0.8) & 0.174 (0.0207) & 0.172 (0.0106) & 0.172 (0.0072) & 0.175 (0.0051) \\
(0.8, 0.5) & 0.391 (0.0131) & 0.418 (0.0124) & 0.399 (0.0075) & 0.404 (0.0042) \\
(0.8, 0.8) & 0.650 (0.0344) & 0.654 (0.0186) & 0.635 (0.0151) & 0.654 (0.0079) \\
\bottomrule
\end{tabular}
\end{table}

The results for Algorithm \ref{alg:WP_unknown_sigma} when the scaling parameter is not known is given below. 

\begin{table}[H]
\centering
\caption{Bifractional Brownian motion: Algorithm \ref{alg:WP_unknown_sigma} with unknown $\sigma^2$}
\label{tab:dpw-bfbm-convergence}
\begin{tabular}{c ccccc}
\toprule
$(H,K)$ & $N=128$ & $N=256$ & $N=512$ & $N=1024$ & $N=8192$ \\
\midrule
$(0.2,\,0.5)$ &
0.230 (0.0627) &
0.200 (0.0428) &
0.161 (0.0192) &
0.149 (0.0139) &
0.122 (0.0058) \\
$(0.8,\,0.5)$ &
0.431 (0.0280) &
0.422 (0.0301) &
0.417 (0.0196) &
0.411 (0.0130) &
0.411 (0.0070) \\
$(0.2,\,0.8)$ &
0.266 (0.0558) &
0.235 (0.0332) &
0.212 (0.0200) &
0.199 (0.0112) &
0.186 (0.0081) \\
$(0.8,\,0.8)$ &
0.639 (0.0317) &
0.616 (0.0287) &
0.625 (0.0245) &
0.607 (0.0174) &
0.619 (0.0094) \\
\bottomrule
\end{tabular}
\end{table}
Like before, the data visualization is given. 

\begin{figure}[H]
\centering
\includegraphics[width=0.45\linewidth]{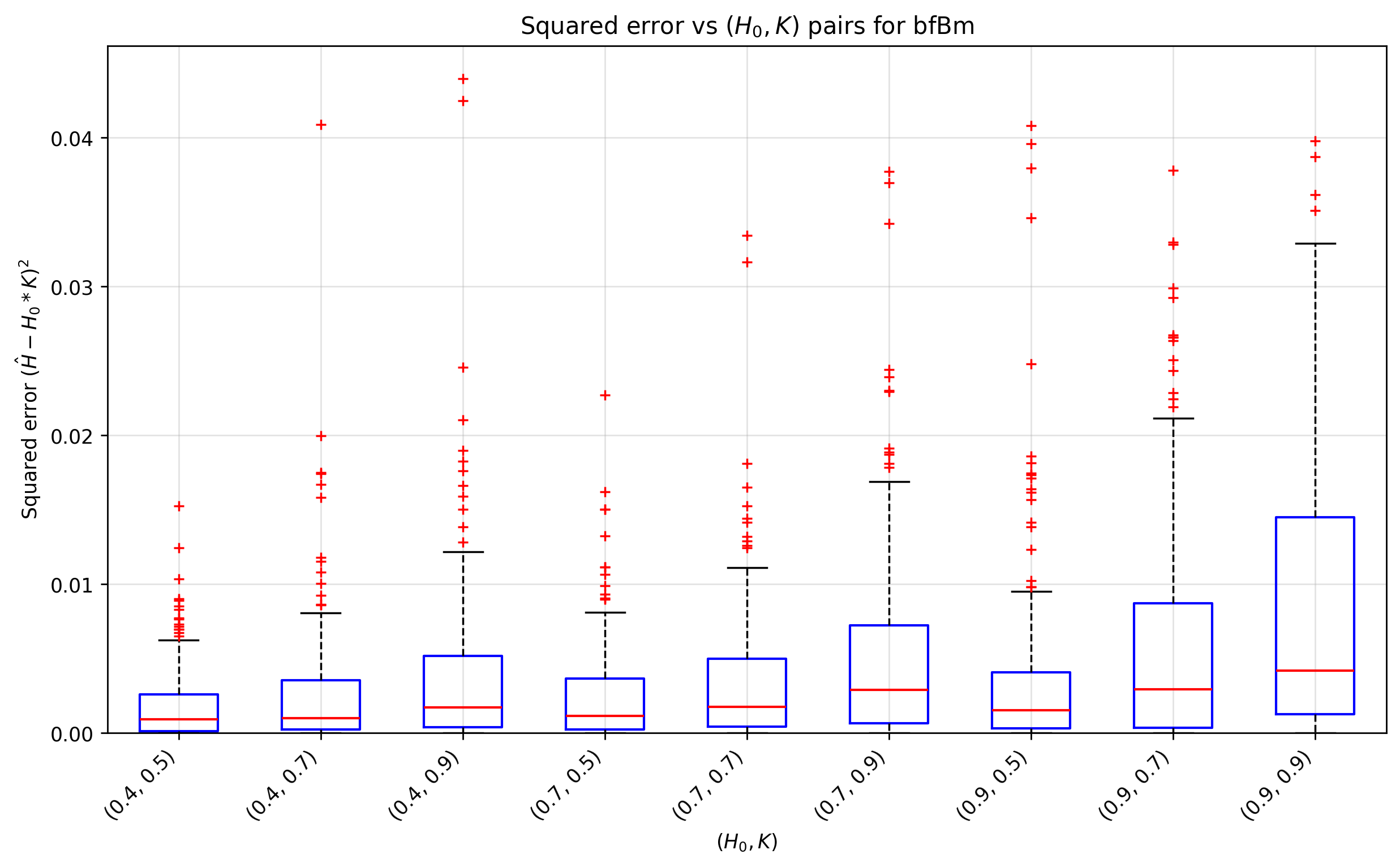}
\includegraphics[width=0.45\linewidth]{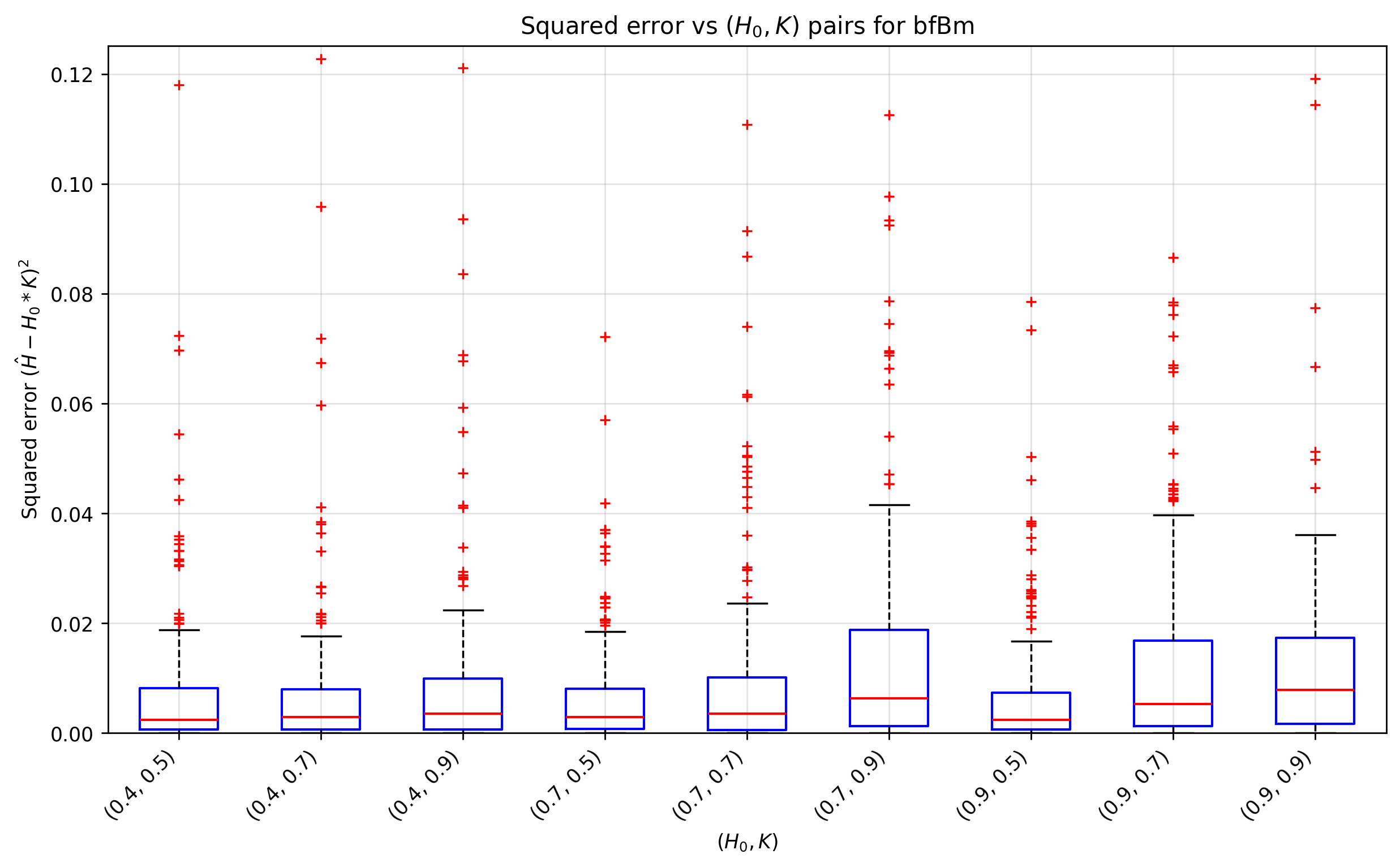}
\caption{The results are given by left: Algorithm \ref{alg:WP_known_sigma}, right: Algorithm \ref{alg:WP_unknown_sigma}. Each $H_0$ consists of 200 simulated paths of length $N=4096$. Squared error of the estimated parameter is computed for each path and drawn into a boxplot. Outliers over 99th percentile are cut out for the sake of clear visuals.}
\label{fig:boxplot_fbm}
\end{figure}

\subsection{Trifractional Brownian Motion}
Another example of Gaussian processes with non-stationary increments is the trifractional Brownian motion  $\{\widetilde B^{H,K}(t)\}_{t\ge0}$  \citep{ma2013schoenberg}. It is a
zero-mean Gaussian process with $\widetilde B^{H,K}(0)=0$ a.s. and covariance function: for $s,t\ge0$,
\begin{equation}
\label{def:tfBm}
\mathbb Cov\left(\widetilde B^{H,K}(t),\widetilde B^{H,K}(s)\right)=\sigma^2\left(t^{2HK}+s^{2HK}-\left(t^{2H}+s^{2H}\right)^K\right),
\end{equation}
where the parameters $H\in(0,1)$ and $K\in(0,1)$. From its definition, we see that $\{\widetilde B^{H,K}(t)\}_{t\ge0}$  is $HK$-self-similar, and its increment processes are not stationary. Like the bifractional Brownian motion, this process also cannot be simulated using the sum of stationary increments. So far as we know, the only simulation approach suggested in the literature is the DPW method.  When $\sigma^2=1$, the estimation of the self-similarity index $HK$ is then obtained by a slight modification of Algorithm \ref{alg:WP_known_sigma}: 
\begin{alg}
 \label{alg:trifBm}
Let $(\widetilde B^{H,K}(j/n))_{j=1,\ldots,n}$ be an observed sample path of $\{\widetilde B^{H,K}(t)\}_{t\in[0,1]}$ with covariance function $\gamma_{\widetilde B^{H,K}}$. 
 \begin{description}

 \item[Step 1] Let the sequences 
  $a=(a_0,\ldots,a_n)
 $ and $b=(b_0,\ldots,b_n)$: 
$$
a=\left(\widetilde B^{H,K}\left(\frac{\lfloor n^{j/n}\rfloor}{n}\right)\right)_{j=0,\ldots,n}~\mbox{and}~b=\left(n^{1-j/n}\right)_{j=0,\ldots,n},
$$
where  $\lfloor\bullet\rfloor$ denotes the floor number. 
  \item[Step 2] Use a non-linear least squares method to solve $H,K$ from the equations
$$
f_n(H,K):=\frac{\sum_{j=0}^{n}a_j^2b_j^{2HK}}{n+1}-(2-2^K)=0.
$$
\end{description}
\end{alg}

Like bifractional Brownian motion, the same parameters for testing are applied to tfBm simulated by the DPW algorithm. The results are given below:

\begin{table}[H]
\centering
\caption{Trifractional Brownian motion: Algorithm \ref{alg:trifBm}.}
\label{tab:tfbm}
\begin{tabular}{c cccc}
\toprule
$(H,K)$ & $N=128$ & $N=256$ & $N=512$ & $N=1024$ \\
\midrule
(0.2, 0.5) & 0.203 (0.1038) & 0.161 (0.0700) & 0.170 (0.0549) & 0.171 (0.0442) \\
(0.2, 0.8) & 0.263 (0.0904) & 0.263 (0.0603) & 0.212 (0.0555) & 0.212 (0.0312) \\
(0.8, 0.5) & 0.442 (0.0672) & 0.461 (0.0460) & 0.427 (0.0319) & 0.448 (0.0195) \\
(0.8, 0.8) & 0.659 (0.0520) & 0.655 (0.0327) & 0.549 (0.0293) & 0.656 (0.0188) \\
\bottomrule
\end{tabular}
\end{table}

Finally the results for trifractional Brownian motion when the scaling parameter is not known, given by Algorithm \ref{alg:WP_unknown_sigma}, is given below.

\begin{table}[H]
\centering
\caption{Trifractional Brownian motion: Algorithm \ref{alg:WP_unknown_sigma}}
\label{tab:dpw-tfbm-convergence}
\begin{tabular}{c ccccc}
\toprule
$(H,K)$ & $N=128$ & $N=256$ & $N=512$ & $N=1024$ & $N=8192$ \\
\midrule
$(0.2,\,0.5)$ &
0.134 (0.0194) &
0.140 (0.0210) &
0.136 (0.0236) &
0.130 (0.0163) &
0.121 (0.0094) \\
$(0.8,\,0.5)$ &
0.384 (0.0414) &
0.421 (0.0375) &
0.408 (0.0399) &
0.406 (0.0303) &
0.405 (0.0185) \\
$(0.2,\,0.8)$ &
0.212 (0.0367) &
0.182 (0.0220) &
0.193 (0.0218) &
0.179 (0.0173) &
0.172 (0.0144) \\
$(0.8,\,0.8)$ &
0.616 (0.0434) &
0.626 (0.0379) &
0.628 (0.0412) &
0.642 (0.0299) &
0.628 (0.0194) \\
\bottomrule
\end{tabular}
\end{table}

The error analysis result is illustrated below:

\begin{figure}[H]
\centering
\includegraphics[width=0.45\linewidth]{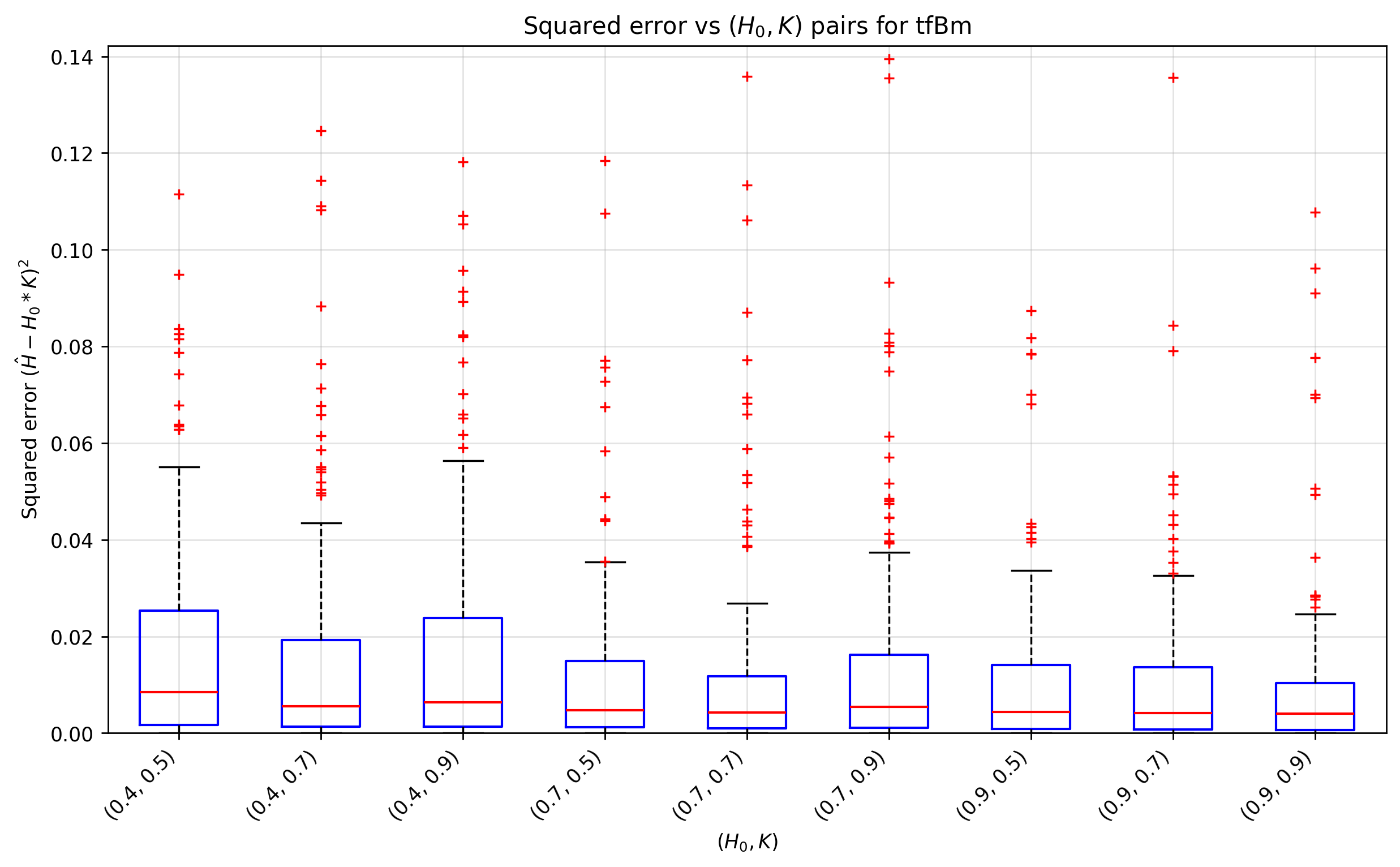}
\includegraphics[width=0.45\linewidth]{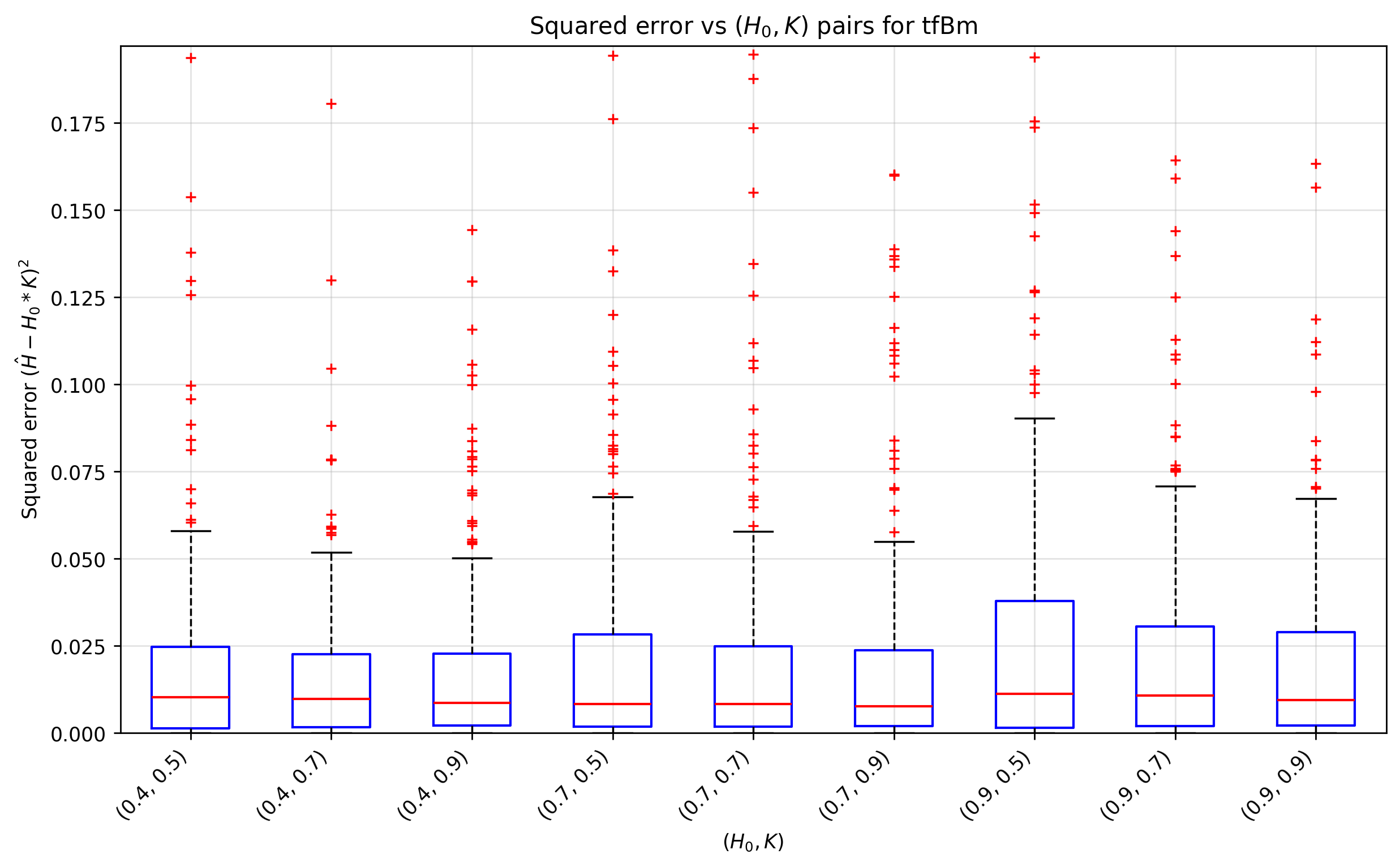}
\caption{The results are given by left: Algorithm \ref{alg:WP_known_sigma}, right: Algorithm \ref{alg:WP_unknown_sigma}. Each $H_0$ consists of 200 simulated paths of length $N=4096$. Squared error of the estimated parameter is computed for each path and drawn into a boxplot. Outliers over 99th percentile are cut out for the sake of clear visuals.}
\label{fig:boxplot_fbm}
\end{figure}
We conclude that:
\begin{description}
\item[(1)] Our estimation method is competitive to the literature ones for estimation of the Hurst parameter of fractional Brownian motion.
\item[(2)] Our estimation fills the gaps when the target process has non-stationary increments. This makes our method be able to estimate the self-similarity indexes of sub-fBm, bifBm, and trifBm.
\end{description}
\subsection{Code availability}
In this section, we list the codes for all methods mentioned in this paper. Among them, we have contributed the following:
\begin{itemize}
    \item Estimation process to obtain comparison data (Python): \url{https://github.com/william11074/Self-similar-Statistical-Inference}
    \item Simulation algorithm for fBm, sfBm, bfBm, and tfBm using both Wood Chan's method and DPW method, along with implementation of the QV method: pip install fractal-analysis:\\ \url{https://pypi.org/project/fractal-analysis/}
\end{itemize}

\section{Real World Application: The Nile River}
Studies conducted by Bardet \cite{bardet2000testing} and Lee et. al \cite{lee2016testing} show that the annual minimum water levels of the Nile River are self-similar with the property of long-range dependency. These papers detail that the self-similarity index of the observed water levels between years 722 and 1281 match $H\approx0.88$. \\

We take the interval from the year 900 to 1200 for the Nile River data and using Algorithm \ref{alg:WP_unknown_sigma}, we get a predicted $H$ value of $0.8545$. However, applying parametric bootstrap bias correction, the corrected value $H_{BC}=0.8790$, matching the results in literature. Open source data is used from \url{http://lib.stat.cmu.edu/S/beran}.

\begin{figure}[H]
\centering
\includegraphics[width=0.8\linewidth]{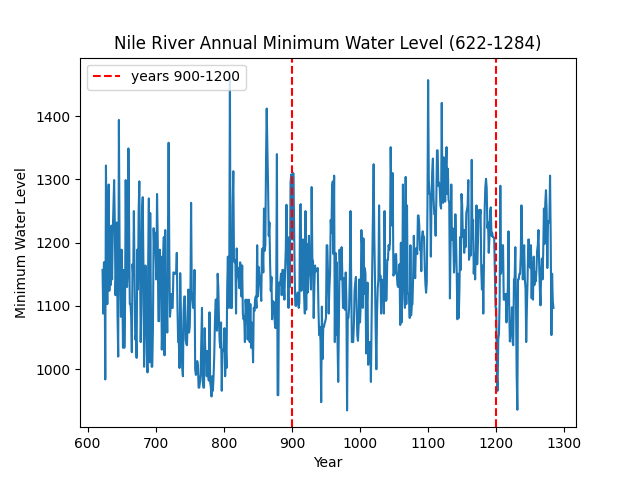}
\caption{Plotted is the available dataset of annual Nile River data from years 622-1284 C.E. The years 900-1200 is chosen for the estimation of the self-similarity index.}
\end{figure}

\section{Conclusion}
This paper develops a general way to estimate the self-similarity index for a given second-order self-similar stochastic process that is ergodic with emphasis that the method works for non-stationary increments self-similar processes. Two estimation algorithms are given for when the scaling parameter is known or unknown. Simulation and testing show that results are accurate and comparable to other methods in the literature. In addition, application to real-world data, the NIle River's historical water level, shows use cases. GitHub code is documented and provided for the estimation algorithms. 

\bibliographystyle{imsart-number.bst}
\bibliography{ml}

\end{document}